\documentclass{amsart}
\usepackage{amssymb}
\usepackage{amsmath}
\usepackage{amscd}
\usepackage{graphicx}

\theoremstyle{plain}
\newtheorem{theorem}{Theorem}[section]
\newtheorem{lemma}[theorem]{Lemma}
\newtheorem{corollary}[theorem]{Corollary}
\newtheorem{proposition}[theorem]{Proposition}
\newtheorem{observation}[theorem]{Observation}

\theoremstyle{definition}
\newtheorem{definition}[theorem]{Definition}

\newtheorem{example}[theorem]{Example}
\newtheorem{notation}[theorem]{Notation}
\newtheorem{standing assumption}[theorem]{Standing Assumption}
\newtheorem*{st assumption}{Standing Assumption}

\theoremstyle{remark}
\newtheorem{remark}[theorem]{Remark}
\newtheorem*{ack}{Acknowledgment}

\begin{document}
\title[Invariants for one-dimensional spaces]
{Ordered group invariants for\\
one-dimensional spaces}

\author{Yi, Inhyeop}

\address{Department of Mathematics, 
University of Maryland, 
College park, MD, 20742-4015, USA}

\email{inhyeop{@}math.umd.edu}



\begin{abstract}
We show that the Bruschlinsky group with
the winding order is a homeomorphism invariant
for a class of one-dimensional inverse limit spaces.
In particular
we show that
if a presentation of
an inverse limit space satisfies
the Simplicity Condition, then
the Bruschlinsky group with the winding order of
the inverse limit space is a dimension group
and is a quotient of 
the dimension group with the standard order
of the adjacency matrices
associated with the presentation.
\end{abstract}
\maketitle

\section{Introduction}\label{s1} 
Ordered groups have been useful invariants for 
the classification of many different categories.
A class of ordered groups, {\bf dimension groups},
was used in the study of
$C^*$-algebras to classify AF-algebras (\cite{ef}),
and Giordano, Herman, Putnam and Skau (\cite{gps, hps})
defined (simple) dimension groups 
in terms of dynamical concepts
to give complete information about
the orbit structure of
zero-dimensional minimal dynamical systems.
Swanson and Volkmer (\cite{sv}) showed that
the dimension group of a primitive matrix is
a complete invariant for {\it weak equivalence},
which is called {\it $C^*$-equivalence}
by Bratteli, J{\o}rgensen, Kim, and Roush (\cite{bjkr}). 
And Barge, Jacklitch, and Vago (\cite{bjv}) showed that,
for certain class of one-dimensional
inverse limit spaces,
two spaces are homeomorphic if and only if
their associated substitutions are weak equivalent,
and if two inverse limit spaces are homeomorphic
and the square of their connection maps are
orientation preserving,
then the dimension groups
of the adjacency matrices associated to
the substitutions are order isomorphic.

A recent development (\cite{bd,bjv,bh,fo,gps,sv}) is
the refinement of $\check H^1(X)$ as a topological invariant
for certain one-dimensional spaces $X$,
by making this group an ordered group.
Here $\check H^1(X)$ is the direct limit of first
cohomology groups on graphs approximating the space X.
There is a natural order on the first cohomology of a
graph (a coset is positive if it contains a nonnegative
function), and the {\it standard order} on $\check H^1(X)$
is the direct limit order derived from the natural
graph orders (see definition \ref{3.remark}).
Except for parts of \cite{bh} and \cite{fo},
the ordered cohomology results
have involved the standard order.

A second order on $\check H^1(X)$, the 
{\it winding order}, 
is geometrically natural as its 
positive elements are the homotopy classes of
continuous orientation preserving maps
from $X$ to $S^1$. 
Boyle and Handelman (\cite{bh})
defined the {\bf winding order} for 
suspension spaces of zero-dimensional dynamical systems, 
and showed that in some (but not all) cases it agrees 
with the standard order. 
Forrest (\cite{fo}) defined the winding order
for the first {\v C}ech cohomology groups of 
directed graphs (thus taking the step of removing 
dynamics), and used this 
 \cite{fo}
to show that 
whenever  two one-dimensional inverse limit spaces
 are  {\it pro-homotopy equivalent},
then their first {\v C}ech cohomology groups
with the standard order are order isomorphic.

In this paper, we extend the definition of 
the winding order to a 
large class of one dimensional spaces, 
\lq\lq compact branched matchbox manifolds\rq\rq .  
We show that for a compact connected
orientable branched matchbox manifold with 
an inverse limit presentation satisfying 
 the Simplicity Condition, 
the Bruschlinsky group with the winding order is
a simple dimension group, and the winding order 
equals the standard order.  
This is a natural extension
of the relations between zero-dimensional
minimal systems and simple dimension groups in
Giordano, Herman, Putnam and Skau (\cite{gps, hps})
to an appropriate class of one-dimensional 
spaces. 
As a corollary
we obtain an independent proof of some results of
Forrest and Barge, Jacklitch and Vago (\cite{bjv, fo}) 
computing dimension group invariants for the 
oriented generalized one-dimensional solenoids of Williams 
(\cite{w1,w2,yi}).

The outline of the paper is as follows.
In section 2, using works of
Aarts and Oversteegen (\cite{ao}),
Marde{\v s}i{\' c} and Segal (\cite{ms})
and Rogers (\cite{ro}),
we define  
compact connected orientable
branched matchbox manifolds, and show they all have 
presentations by orientation preserving maps
of finite directed nondegenerate graphs.
In section 3,
we show that the Bruschlinsky group
with the winding order of
a compact connected orientable
branched matchbox manifold with
the Simplicity Condition
is order isomorphic to the direct limit of
the graph groups with the standard order
defined from the presentation (and therefore 
the winding and standard orders agree). 
And in section 4, we recall the axioms
for one-dimensional generalized solenoids
and calculate the Bruschlinsky groups with
the winding order of an example in which 
 the Bruschlinsky group is not given by 
the obvious direct limit of presenting matrices.

\section{Branched matchbox manifolds
and ordered groups}\label{s2}
Aarts and Oversteegen (\cite{ao})
defined a {\it matchbox manifold} to be 
a separable metric space $Y$ such that
each point $y\in Y$ has 
a neighborhood which is homeomorphic to
$S_y\times I_y$
where $S_y$ is a zero-dimensional space and
$I_y$ is an open interval.
For a topological embedding
$h\colon S_y\times I_y\to Y$,
they called $h\left(S_y\times I_y\right)$
a {\it matchbox neighborhood} of $y\in Y$.
A matchbox manifold $Y$ is called {\it orientable}
if each arc component $C_{\alpha}$, $\alpha\in A$,
of $Y$ has a parameterized immersed arc
$p_{\alpha}\colon {\mathbb R} \to C_{\alpha}$
such that each point $y\in Y$
has a matchbox neighborhood $h\left(S_y\times I_y\right)$
with the following property:
for each $\alpha \in A$ and each $t\in \mathbb R$
with $p_{\alpha}(t)\in h\left(S_y\times I_y\right)$
there exists an open interval $I$ containing $t$
such that $pr_2\circ h^{-1}\circ p_{\alpha}$
is increasing on $I$
where $pr_2$ is the canonical projection 
from $S_{p_{\alpha}(t)}\times I_{p_{\alpha}(t)}$
to $I_{p_{\alpha}(t)}$.

\begin{theorem}[{\cite{ao}}]\label{2.2}
For a one-dimensional space $Y$,
the following are equivalent:
\begin{enumerate}
\item[(1)]
$Y$ is an orientable matchbox manifold.
\item[(2)]
$Y$ is the phase space of a flow without rest point.
\item[(3)]
There exists a cross section $K$
with return time map $r_K$ such that
$Y$ is the standard suspension of $(K,r_K)$.
\end{enumerate}
\end{theorem}

\subsection*{Branched matchbox manifold}
We define a {\it branched matchbox} to be
a topological space homeomorphic to
$U =\bigl(\left(S_1\times \left(-1,0\right]\right)
\cup \left(S_2\times \left[0,1\right)\right)\bigr)/\sim$ 
such that $S_1$ and $S_2$ are
zero dimensional separable metrizable spaces
and there is a (closed) equivalence 
relation $\approx$ on $S_1 \cup S_2$ such that 
\begin{enumerate}
\item[(1)]
For every $s_1\in S_1$
($\sigma_2\in S_2$, respectively)
there exists at least one $s_2\in S_2$
($\sigma_1\in S_1$, respectively)
such that $s_1\approx s_2$
($\sigma_1\approx \sigma_2$, respectively),
\item[(2)]
$(S_1 \cup S_2)/{\approx}$ is
a zero dimensional metrizable space
with the quotient topology and
\item[(3)]
$(s_1,i)\sim (s_2,j)$ if and only if
either $s_1 \approx s_2$ and $i=j=0$
or $s_1=s_2$ and $i=j$.
\end{enumerate}
\begin{remark}
In this paper,
we will always be concerned with the case that 
$S_1$ and $S_2$ are compact.
\end{remark}

%

For $s_1\in S_1$ and $s_2\in S_2$ such that
$s_1 \approx s_2$,
the set
$$
\bigl(\left(\{s_1\}\times \left(-1,0\right]\right)
\cup \left(\{s_2\}\times \left[0,1\right)\right)\bigr)/\sim
$$ 
is called a {\it match}.

A 
{\it branched matchbox manifold}
is a separable metrizable space $Y$
together with a collection of maps
called {\it charts}
such that
\begin{itemize}
\item[(1)]
a chart is a homeomorphism
$h\colon V\to U$
where $V$ is an open set in $X$ and
$U$ is a branched matchbox,
\item[(2)]
every point in $Y$ is in the domain of some chart,
\item[(3)]
for charts $h_1\colon V_1\to U_1$ and
$h_2\colon V_2\to U_2$
the change of coordinates map
$h_2\circ h_1^{-1}\colon
h_1\left(V_1\cap V_2\right) \to  
h_2\left(V_1\cap V_2\right)$
is continuous. 
\end{itemize}

Every branched matchbox $U$ has the direction given
by the second coordinate with a continuous projection
$p_U\colon U\to (-1,1)$ defined by
$[(z,j)]\mapsto j$.
Following the approach of Aarts and Oversteegen
(\cite[\S3]{ao}),
we call a branched matchbox manifold $Y$
{\it orientable}
if it can be covered by branched matchboxes
with directions agreeing on overlaps, i.e.,
there are oriented branched matchboxes
$U_i$ with projections $p_i\colon U_i \to (-1,1)$,
open sets $V_i$ covering  $Y$,
and homeomorphisms $h_i\colon V_i\to U_i$ such that
for every $i,j$ and every locally one-to-one curve
$\gamma \colon  [0,1]\to V_i \cap V_j$,
$p_i\circ h_i\circ \gamma$ is increasing
if and only if
$p_j\circ h_j\circ \gamma$ is increasing. 
The particular collection of charts,
maximal with respect to
this change of coordinate property,
is called an {\it orientation}
of the branched matchbox manifold $Y$.

\subsection*{Ordered group}
A {\it preordered group} is a pair $(G,G_+)$
where $G$ is an Abelian group
and the {\it positive cone}
$G_+$ is a submonoid of $G$ which generates $G$.
We write $g_1\le g_2$ if $g_2-g_1\in G_+$ for 
$g_1,g_2\in G$.
If $(G,G_+)$ satisfies the additional condition that
$G_+\cap -G_+=\{ 0 \}$,
then  $(G,G_+)$ is called an {\it ordered group}.

An {\it order unit} in a preordered group
is an element $u\in G_+$ such that
for every $g\in G$ there exists a positive integer
$n=n(g)$ such that $g\le nu$.
A preordered group $(G,G_+)$ is {\it unperforated}
if for every $g\in G$ and positive integer $n$,
$ng\in G_+$ implies $g\in G_+$.
We say that an ordered group $(G,G_+)$
satisfies the {\it Riesz Interpolation property}
if given $g_1,g_2, h_1,h_2 \in G$ with
$g_i\leq h_j$ ($i,j=1,2$),
then there is a $k\in G$ such that
$g_i\leq k \leq h_j$.

\subsection*{Bruschlinsky group with the winding order}
For a compact metric space $Y$,
let $C(Y,S^1)$ be the set of continuous functions
from $Y$ to $S^1$, and
$$
R(Y)=\{ \phi\in C(Y,S^1)\mid \phi(y)=\exp(2\pi i g(y)) 
\text{ for some } g\in C(Y,\mathbb{R})\}.
$$  
Then $R(Y)$ is 
the subgroup of functions
homotopic to a constant map in $C(Y,S^1)$.
The {\it Bruschlinsky group} of $Y$ (\cite[\S4.3]{pt})
is given by
$$
\text{Br}(Y)=C(Y,S^1)/R(Y).
$$
It is well-known  that ${\check H}^1(Y)$,
the first {\v C}ech cohomology group of $Y$, 
is isomorphic to the Bruschlinsky group of $Y$
(\cite{bh, kr}).

Now suppose that $Y$ is
an oriented compact branched matchbox manifold.
Let $C_{\oplus}(Y,S^1)$ be the set of
$\phi\in C(Y,S^1)$ such that there exists a map
$\psi \in R(Y)$ such that
$\phi \cdot \psi$ is non-orientation-reversing, i.e.,
for every orientation preserving parameterized curve
$\gamma\colon \mathbb{R}\to Y$,
$(\phi \cdot \psi)(\gamma)(t)$
does not move in the clockwise direction
as $t\in \mathbb{R}$ increases.

Define 
$ \text{Br}_{\oplus}(Y)
=\bigl\{ [\phi]\mid \phi\in C_{\oplus}(Y,S^1)\bigr\}$.
Then $(\text{Br}(Y),\text{Br}_{\oplus}(Y))$ is 
a preordered group.
We call this preorder
the {\bf winding order} (\cite[\S4]{bh}).

\begin{remark}[{\cite[4.7]{bh}}]
It is possible that
the Bruschlinsky group with the winding order
of compact orientable space is not
an ordered group.
\end{remark}

\begin{observation}\label{1.1}
Homeomorphic orientable compact metric spaces
have order-isomorphic Bruschlinsky groups
with the winding order.
\end{observation}

\begin{proposition}[{\cite{kr}}]\label{2.k}
The Bruschlinsky group of a 
compact branched matchbox manifold
is a torsion-free group.
\end{proposition}

Recall that a continuum is
a compact connected metric space.

\begin{lemma}[{\cite{kr}}]\label{2.l}
Let $Y$ be a continuum, $\phi\in C(Y,S^1)$,
and $p_n\colon S^1 \to S^1$ defined by $z\mapsto z^n$
for every positive integer $n$.
Then $n\cdot \left[\phi\right]=\left[ p_n\circ \phi\right]$.
\end{lemma}

\begin{proposition}\label{2.r}
The Bruschlinsky group with the winding order of
a compact connected oriented branched matchbox manifold $Y$
is unperforated.
\end{proposition}
\begin{proof}
Suppose that $\phi\in C(Y,S^1)$
and $n\in \mathbb{Z}_+$ such that 
$n\cdot \left[\phi\right] =\left[p_n\circ \phi \right]
\in \text{Br}_{\oplus}(Y)$.
Then there exists a map
$\psi\in R(Y)$ given by
$y\mapsto \exp\left( 2\pi i g(y)\right)$
with $g\in C(Y,\mathbb{R})$ such that
$\left( p_n\circ \phi \right) \cdot \psi$
is non-orientation reversing.

Define $\tilde{\psi}\colon Y\to S^1$
by $y\mapsto \exp\left( 2\pi i \cdot \frac{1}{n}g(y)\right)$.
Then we have $\tilde{\psi}\in R(Y)$ and
$\left( p_n\circ \phi \right) \cdot \psi
=p_n\circ \left( \phi \cdot \tilde{\psi}\right)$.
For every orientation preserving parameterized curve
$\gamma\colon \mathbb{R}\to Y$,
$$
\left(\left( p_n\circ \phi \right) \cdot \psi\right)
\circ \gamma(t)
=p_n\circ
\left( \phi \cdot \tilde{\psi}\right)\circ \gamma(t)
=p_n\circ\left(
\left( \phi \cdot \tilde{\psi}\right)\circ \gamma(t)
\right)
$$
does not move clockwise on $S^1$
as $t\in \mathbb{R}$ increases.
So $\phi \cdot \tilde{\psi}$
is non-orientation-reversing as
$n$ is a positive integer.
Therefore
$\phi\in C_{\oplus}(Y,S^1)$,
and $\left(\text{Br}(Y),\text{Br}_{\oplus}(Y)\right)$
is unperforated.
\end{proof}

\begin{remark}
If $Y$ is a compact connected orientable matchbox manifold,
then the above Propositions \ref{2.k} and \ref{2.r}
follow from Propositions 4.5 and 3.4 of \cite{bh}
and Theorem \ref{2.2}.
\end{remark}

\subsection*{One dimensional continua}
In \cite{ro}, Rogers introduced
the following notations for
one-dimensional continua.

Suppose that $X_1$ and $X_2$ are graphs and 
that $\mathcal{V}_i$ and $\mathcal{E}_i$ are
the vertex set and the edge set
of $X_i$, respectively, $i=1,2$.
A continuous onto map $f\colon X_2\to X_1$ is called 
{\it simplicial relative to} $(\mathcal{V}_1,\mathcal{V}_2)$
if $f(\mathcal{V}_2)\subseteq \mathcal{V}_1$ and
for every edge $e_2\in \mathcal{E}_2$
there is an edge $e_1\in \mathcal{E}_1$ such that
$f|_{e_2\backslash \mathcal{V}_2}$
is a homeomorphism onto ${e_1\backslash \mathcal{V}_1}$
or a constant map.
The map $f\colon X_2\to X_1$ is {\it simplicial}
if it is simplicial relative to some vertex sets of
$X_1$ and $X_2$. 
And $f$ is called {\it light} 
if the preimage of each point is totally disconnected.

An inverse limit sequence $\{X_k,f_k\}$ on graphs is
called {\it light simplicial}
if each $f_k$ is light simplicial,
and is called {\it light uniformly simplicial} if
each $X_k$ is a graph with a vertex set $\mathcal{V}_k$
and each map $f_k\colon X_k\to X_{k-1}$ is light simplicial
relative to $(\mathcal{V}_{k-1},\mathcal{V}_k)$.

\begin{theorem}[{\cite{ms,ro}}]\label{2.3} 
Suppose that $\overline{X}$ is
a one-dimensional continuum.
\begin{enumerate}
\item[(1)]
$\overline{X}$ is homeomorphic to an inverse limit of
a light simplicial sequence $\{X_k,f_k\}$ on graphs, and
\item[(2)]
$\overline{X}$ is homeomorphic to
a light uniformly simplicial inverse limit on graphs 
if and only if there exists a map
$\pi:\overline{X}\to [0,1]$ such that
$\pi^{-1}(\{0,1\})$ is totally disconnected
and $\pi|_{e}$ is a homeomorphism for every $e$
which is a closure of component of
$\overline{X}\backslash \pi^{-1}(\{0,1\})$.
\end{enumerate}
\end{theorem}

Suppose that $\{ X_k,f_k\}$
is a light simplicial sequence on graphs.
Let 
$$\overline{X}
=X_0\overset{f_1}{\longleftarrow}
X_1\overset{f_2}{\longleftarrow}\cdots
=\bigl\{(x_0,x_1,x_2,\dots )\in \prod_0^\infty X_k\, 
                   |\, f_{k+1}(x_{k+1})=x_k \bigr\}.
$$
For a one-dimensional continuum $Y$,
we call the sequence $\{ X_k,f_k \}$
a {\bf presentation} of $Y$
if $\overline{X}$ is homeomorphic to $Y$.

\begin{notation}
Suppose that $G$ is a directed graph.
We consider a directed edge $e$ of $G$
as the image of a local homeomorphism
from $[0,1]$ to $e$ such that
$e(0)$ is the initial point of $e$
and $e(1)$ is the terminal point.
Then we can represent each point
$x\in e$ as $e(t)$ (possibly $e(0)=e(1)$).
\end{notation}

Recall that a continuous map $p\colon [0,1]\to G$,
a directed graph,
is {\it orientation preserving}
if $e^{-1}\circ p\colon I\to [0,1]$ is increasing
for every interval $I\subset [0,1]$
such that $p(I)$ is a subset of a directed edge $e$.
A continuous map $f\colon G_1\to G_2$
between two directed graphs
is {\it orientation preserving}
if, for every orientation preserving map
$p \colon [0,1]\to G_1$,
$f\circ p \colon [0,1]\to G_2$ is
orientation preserving (\cite{fo}). 
A directed graph is called {\it nondegenerate}
if every vertex has at least one incoming edge and
at least one outgoing edge.

Suppose that $Y$ is
a compact connected oriented branched matchbox manifold.
Since $Y$ is a one-dimensional continuum,
there is a light  simplicial presentation
$\{X_k, f_k\}$ of $Y$ by Theorem \ref{2.3}.
The following proposition shows that
the orientation of $Y$ decides
the directions of edges in each coordinate space
$X_k$ so that every connection map
$f_k\colon X_k\to X_{k-1}$
is orientation preserving.

\begin{proposition}\label{2.4}
Suppose that $Y$ is
a compact connected oriented branched matchbox manifold.
Then $Y$ has a light  simplicial presentation 
by orientation preserving maps of
directed nondegenerate graphs.
\end{proposition}
\begin{proof}
Suppose that 
$\{ h_U\colon V\to U\}$ is an orientation of $Y$
where $U$ is a branched matchbox
with the projections $p_U\colon U \to (-1,1)$.
Let $\{X_k, f_k\}$ be
a light uniformly simplicial presentation
of $Y$ given by Theorem \ref{2.3},
and $\pi_k\colon Y \to X_k$  the canonical
projection to the $k$th coordinate space.
If $e$ is an edge of $X_k$ with
$\pi_k^{-1}\left(e\backslash \mathcal{V}_k\right)
\cap h^{-1}_U(U)\ne \emptyset$,
then give the direction to the set
$\left(e\backslash \mathcal{V}_k\right)\cap
\left(\pi_k\circ h^{-1}_U(U)\right)\subset e$
so that, 
for every curve $\gamma\colon [0,1] \to
\pi_k^{-1}\left(e\backslash \mathcal{V}_k\right)
\cap h^{-1}_U(U)$,
$p_U\circ h_U\circ\gamma$ is increasing
if and only if
${e}^{-1}\circ \pi_k\circ \gamma$ is increasing.
Since $\{h_U\}$ is an orientation of $Y$,
we can extend this direction on
$\left(e\backslash \mathcal{V}_k\right)\cap
\pi_k\circ h^{-1}_U(U)$ to $e$,
and each edge $X_k$ has a direction
induced by the orientation of $Y$.

Suppose that
$x=(x_0, x_1,\dots)$ is a point in $Y$ such that
$x_k\in X_k$ is a vertex
and that $U$ is a branched matchbox such that
the domain of $h_U$  contains $x$.
Then there is a match $M\subset U$
containing $h_U(x)$ such that 
$p_U|_M\circ h_U(x)=t$ for some  $t\in (-1,1)$.
Since
$\pi_k\circ h^{-1}_U\circ
\left(p_U|_M\right)^{-1}\left( (-1,t)\right)$
and $\pi_k\circ h^{-1}_U\circ
\left(p_U|_M\right)^{-1}\left( (t,1) \right)$
are nonempty sets in $X_k$,
there exist an edge $e_-$ such that
$\bigl(\pi_k\circ h^{-1}_U\circ
\left(p_U|_M\right)^{-1}\left( (-1,t)\right)\bigr)
\cap e_-\ne \emptyset$,
which is incoming to $x_k$,
and an edge $e_+$ such that
$\bigl(\pi_k\circ h^{-1}_U\circ
\left(p_U|_M\right)^{-1}\left( (t,1) \right)\bigr)
\cap e_+ \ne \emptyset$,
which is outgoing from $x_k$.
Therefore $X_k$ is nondegenerate.

Suppose that
$e_k\in \mathcal{E}_k$ and $e_{k-1}\in \mathcal{E}_{k-1}$
are two edges such that $e_{k-1}=f_k(e_k)$,
and that $h_U\colon V\to U$ is a chart such that
$W=\pi_k\circ h^{-1}_U(U)\cap 
\left( e_k\backslash \mathcal{V}_k\right)\ne
\emptyset$.
Then we have
$f_k\left( W \right)\subset
\pi_{k-1}\circ h^{-1}_U(U)\cap
\left(e_{k-1}\backslash\mathcal{V}_{k-1}\right)$,
and for every every curve
$\gamma\colon [0,1]\to
h^{-1}_U(U)\cap \pi_k^{-1}\left( e_k \backslash \mathcal{V}_k
\right)$,
${e}_k^{-1}\circ \pi_k\circ \gamma
\text{ is increasing}
\iff
p_U\circ h_u\circ\gamma \text{ is increasing}
\iff
{e}_{k-1}^{-1}\circ \pi_{k-1}\circ \gamma 
\text{ is increasing}$.

Let $\gamma\colon [a,b]\to
h^{-1}_U(U)\cap \pi_k^{-1}\left(
e_k \backslash \mathcal{V}_k \right)$
be given by
$\pi_k\circ \gamma(t)={e}_k(t)$.
Then we have
$\pi_{k-1}\circ \gamma(t)=f_k\circ {e}_k(t)$,
and
${e}_{k-1}^{-1}\circ \pi_{k-1}\circ \gamma (t)
={e}_{k-1}^{-1}\circ f_k\circ {e}_k(t)$
is increasing as $t$ is increasing.
Therefore $f_k\colon X_k\to X_{k-1}$
is orientation preserving.
\end{proof}

\begin{corollary}
Suppose that
$Y$ is a compact connected orientable branched matchbox manifold.
Then there is a continuous map
$\pi\colon Y\to S^1$ such that
$\pi^{-1}(1)$ is totally disconnected and
$\pi|_{\ell}$ is an orientation preserving homeomorphism
for every $\ell$ which is an
arc component of $Y\backslash \pi^{-1}(1)$.
\end{corollary}
\begin{proof}
Define $\pi\colon Y \to S^1$
by $x=(x_0,x_1,\dots)\mapsto \exp(2\pi i t)$
where $t\in [0,1]$ is given by
$x_0={e}(t)\in e\in \mathcal{E}_0$.
Then $\pi$ is well-defined and 
$\pi^{-1}(1)=\{x\in Y\mid x_0\in \mathcal{V}_0\}$
is a zero dimensional set.
Since $\ell$,
an arc component of $Y\backslash \pi^{-1}(1)$,
is given by
$\ell=(e_0\backslash \mathcal{V}_0,
e_1\backslash \mathcal{V}_1,\dots)$
where $e_i\in \mathcal{E}_i$,
$\pi\colon \ell \to S^1$ given by
$x=({e}_0(t),{e}_1(t),\dots)\mapsto
\exp(2\pi i t)$ is 
an orientation preserving homeomorphism.
\end{proof}

We have the following proposition from theorem \ref{2.3}.

\begin{proposition}
Every compact connected orientable
branched matchbox manifold has 
a light uniformly simplicial presentation.
\end{proposition}

\begin{standing assumption}\label{2sa}
From now on,
a graph means a finite directed nondegenerate graph.
\end{standing assumption}

\section{Orientable one dimensional
inverse limit spaces}\label{s3}

In this section
we suppose that  $\overline{X}$ is
a compact connected oriented branched matchbox manifold
with a presentation $\{X_k,f_k \}$ such that
each $X_k$ is a graph with a fixed vertex set $\mathcal{V}_k$
and each map $f_k\colon X_k\to X_{k-1}$ is 
an orientation preserving map such that
$f_k\left(\mathcal{V}_k\right)\subset \mathcal{V}_{k-1}$
and $f_k|_{X_k\backslash\mathcal{V}_k}$
is locally one-to-one.
Let $\mathcal{E}_k$ be the set of directed edges
in $X_k$ defined by $\mathcal{V}_k$,
$C(\mathcal{E}_k,\mathbb{Z})$ 
the set of integer-valued functions on $\mathcal{E}_k$,
and $C_+(\mathcal{E}_k,\mathbb{Z})$
the subset of $C(\mathcal{E}_k,\mathbb{Z})$
with range in the nonnegative integers $\mathbb{Z}_+$.
For each vertex $p_i$ of $X_k$,
define the {\it vertex function}
$v_i\in C(\mathcal{E}_k,\mathbb{Z})$ 
such that for every edge $e\in \mathcal{E}_k$
\begin{equation*}
v_i(e)=\begin{cases} 1& \text{
if $e$ is an edge from $p_i$ to other vertex point},\\
-1 &\text{
if $e$ is an edge from other vertex point to $p_i$},\\
0 &\text{
if $p_i$ is the initial and terminal point of $e$, or
$p_i\notin e$}.
\end{cases}
\end{equation*}
Denote $V_k$ as the set of integral combinations of 
$\{v_i\}\subset C(\mathcal{E}_k,\mathbb{Z})$,
and call an element of $V_k$ a {\it vertex coboundary}.
Define
$$
\mathcal{G}^{k}=C(\mathcal{E}_k,\mathbb{Z})/V_k
\text{ and }
\mathcal{G}^{k}_+=C_+(\mathcal{E}_k,\mathbb{Z})/V_k.
$$
Then
$(\mathcal{G}^{k},\mathcal{G}^{k}_+,\mathbf{1})$
is a unital preordered group.

\begin{notation}
By a {\bf path} in a graph $X$
we mean a finite sequence
$e_1^{s(1)}\cdots e_n^{s(n)}$ of edges such that,
for $1\le i<n$,
$s(i)=\pm 1$ represents the direction of $e_i$
and 
the terminal vertex of $e_i^{s(i)}$
is the initial vertex of $e_{i+1}^{s(i+1)}$.
We write $e^s\in \wp$ if 
$\wp$ is a path and $e$ is an edge such that
$e^s$ is a factor of $\wp$.
A {\bf cycle} is a path $e_1^{s(1)}\cdots e_n^{s(n)}$
such that the terminal vertex of $e_n^{s(n)}$ is
the initial vertex of $e_1^{s(1)}$.
\end{notation}

We say that a function $g$ in $C(\mathcal{E}_k,\mathbb{Z})$
is {\it zero $($nonnegative, {\rm respectively}$)$
on cycles} if
the sum of $g(e)$ over the edges $e$
of every cycle in $X_k$ is zero
(nonnegative, respectively).

\begin{lemma}[{\cite[\S3]{bh}}]\label{3.le1}
Suppose that $g$ is an element of
$C(\mathcal{E}_k,\mathbb{Z})$.
Then
\begin{itemize}
\item[(1)]
$g$ is an element of $V_k$ if and only if 
$g$ is zero on cycles in $X_k$, and
\item[(2)]
$[g]$ is an element of
$C_+(\mathcal{E}_k,\mathbb{Z})/V_k=\mathcal{G}^{k}_+$
if and only if $g$ is nonnegative on cycles.
\end{itemize}
\end{lemma}

Given $g\in C(\mathcal{E}_k,\mathbb{Z})$, 
define a continuous map
$$
\phi_g\colon X_k\to S^1 \text{ by } x\mapsto \exp(2\pi i t g(e))
\text{ for } x=e(t), \; t\in [0,1].
$$
Then $\phi_g$ is well-defined
as every vertex point maps to $1\in S^1$,
and  $\phi_g$ is an element of $C(X_k,S^1)$.

\begin{lemma}\label{3.le2}
Suppose that $g$ is an element of
$C(\mathcal{E}_k,\mathbb{Z})$.
Then 
$g$ is an element of $V_k$ if and only if
$\phi_g$ is homotopic to
a constant function $1$ in $C({X}_k,S^1)$. 
\end{lemma}
\begin{proof}
Suppose that $g$ is an element of $V_k$.
For each vertex function $v_i$ defined
at the vertex $p_i$ of $X_k$,
define a map $h_{sv_i}\colon X_k\to S^1$
for $0\le s \le 1$ by
$$
h_{s v_i}\left( e(t) \right)
=\begin{cases}
e^{2\pi i s t}&
\text{if $e$ is an edge from $p_i$ to another vertex point},\\
e^{-2\pi i s t}&
\text{if $e$ is an edge from another vertex point to $p_i$},\\
e^{2\pi i s}&
\text{if $p_i$ is the initial and terminal point of $e$},\\  
1& \text{otherwise}.
\end{cases}
$$
Then $s\mapsto h_{sv_i}$,
$0\le s \le 1$, is a homotopy
between $\phi_{v_i}$ and $1$.

Now suppose that
$\phi_g$ and $1$ are
homotopic on ${X}_k$.
Since the winding number of the restriction of
$\phi_g$ on every cycle in $X_k$ is
a homotopy invariant and
$\sum\limits_{e\in \ell} g(e)$ is the winding number
for every cycle $\ell$ in $X_k$,
we have that $g$ is zero on every cycle,
and $g$ is an element of $V_k$ by Lemma \ref{3.le1}.
\end{proof}

Therefore we have a well-defined map 
$$
\iota_k\colon
\mathcal{G}^{k}\to \text{Br}({X}_k)
\text{ given by } [g]\mapsto [\phi_g].
$$

\begin{proposition}\label{3.pro}
Let $\iota_k$ be defined as above.
Then $\iota_k$ is an isomorphism of preordered groups
$\left(\mathcal{G}^{k},\mathcal{G}^{k}_+\right)$ and
$\left(\mathrm{Br}(X_k),\mathrm{Br}_{\oplus}(X_k)\right)$.
\end{proposition}
\begin{proof}
Since $\phi_{g+h}=\phi_g\cdot\phi_h$,
$\iota_k$ is a group homomorphism.
By Lemma \ref{3.le2}, 
$\phi_g$ is homotopic to a constant function $1$
if and only if $g$ is a vertex coboundary.
So we have
$\iota_k\colon \mathcal{G}^k\to \text{Br}({X}_k)$ 
is injective.

To obtain an inverse of $\iota_k$,
suppose that $\phi$ belongs to $C(X_k,S^1)$.
Then we can choose a map
$\rho\colon \mathcal{V}_k\to \mathbb{R}$
where $\mathcal{V}_k$ is the vertex set of $X_k$
such that 
$\phi(p)=\phi\left(2\pi i \rho(p)\right)$
for every vertex $p$ of $X_k$.
Define $S_{\rho}\in C(X_k,S^1)$ by
$$
e(t)\mapsto \exp\Bigl(
2\pi i \bigl(
(1-t)\rho(e(0))+t\rho(e(1))
\bigr)
\Bigr),
\text{ } 0 \le t \le 1 .
$$
Then $S_{\rho}$ is homotopic to the constant map $1$
by $H_u=S_{u\rho}$ for $0\le u \le 1$,
$\phi$ is homotopic to $\phi / S_{\rho}$,
and for every vertex $p$ of $X_k$,
$\left( \phi / S_{\rho} \right)(p)=1\in S^1$.

For each edge $e\in\mathcal{E}_k$, 
let $r_{\phi}(e)$ be the number  of times the loop 
$\left( \phi / S_{\rho} \right)(x)$
winds around $S^1$ as $x=e(t)$ moves on $e$.
Since $\left( \phi / S_{\rho} \right)(p)=1\in S^1$
for every vertex $p$ of $X_k$,
$r_{\phi}(e)$ is well-defined for each edge $e$.  
Then $r_{\phi}\colon e\mapsto r_{\phi}(e)$ is
an element of $C(\mathcal{E}_k,\mathbb{Z})$,
and $\phi_{r_{\phi}}$ wraps around $S^1$ the same
number of times as $\phi / S_{\rho}$. 
Therefore $\phi_{r_{\phi}}$ is homotopic to
$\phi / S_{\rho}$,
and $\left[\phi\right] \mapsto [r_{\phi}]$
gives the desired inverse to $\iota_k$.

Clearly if $g\in C(\mathcal{E}_k,\mathbb{Z}_+)$,
then $\left[\iota_k(g)\right]=\left[\phi_g\right]$
is a positive element in the winding order.
Conversely if $\left[\phi_g\right]\in \mathrm{Br}(X_k)$
is a positive in the winding order,
then there exists a map $\psi\in R(X_k)$
such that $\phi_g \cdot \psi$ is
non-orientation-reversing.
It follows that $g$ has to be nonnegative on cycles,
and we have $\left[ g \right]\in \mathcal{G}^k_+$
by lemma \ref{3.le1}.
Therefore $\iota_k$ is an isomorphism of preordered groups.
\end{proof}

Since $f_{k+1}\colon X_{k+1}\to X_k$ is
an orientation preserving map,
if $e$ is an edge in $\mathcal{E}_{k+1}$,
then $f_{k+1}(e)$ is a path $e_1\cdots e_n$ in $X_k$.
Hence  $f_{k+1}$ induces a map 
$$
f^*_{k+1}\colon C(\mathcal{E}_{k},\mathbb{Z})\to
C(\mathcal{E}_{k+1},\mathbb{Z})
\text{ defined by }
g\mapsto g\circ f_{k+1}
$$
where $(g\circ f_{k+1})(e)= \sum\limits_{i=1}^{n} g(e_i)$
such that $f_{k+1}(e)=e_1\cdots e_n$ in $\mathcal{E}_k$.
And $f_{k+1}$ induces another map 
$$
\tilde{f}^*_{k+1}\colon
C({X}_{k},S^1)\to C({X}_{k+1},S^1)
\text{ defined by }
\phi \mapsto \phi\circ f_{k+1}.
$$

\begin{lemma}\label{3.0}
Let $f^*_{k+1}$ and $\tilde{f}^*_{k+1}$
be given as above.
Then there are well-defined homomorphisms
from $\mathcal{G}^k \text{ to } \mathcal{G}^{k+1}$
and from
$\mathrm{Br}(X_k) \text{ to }\mathrm{Br}(X_{k+1})$
defined by $f^*_{k+1}$ and $\tilde{f}^*_{k+1}$,
respectively.
\end{lemma}
\begin{proof}
For every $v\in V_k$ and every cycle $\ell$ in $X_{k+1}$,
$f_{k+1}(\ell)$ is a cycle in $X_{k}$ and
$f_{k+1}^*(v)(\ell)=v\left( f_{k+1}(\ell) \right)=0$
by Lemma \ref{3.le1}.
Therefore $f_{k+1}^*(v)$ is an element of $V_{k+1}$,
and the map $\mathcal{G}^k \to \mathcal{G}^{k+1}$
given by $[g]\mapsto [f^*_{k+1}(g)]$ is
a well-defined homomorphism.
That $\tilde f_{k+1}^*$ induces
a homomorphism follows from
the definition of the Bruschlinsky group.
\end{proof}

Let us denote these well-defined homomorphisms as
$f_{k+1}^*$ and $\tilde{f}^*_{k+1}$, respectively,
if they do not give any confusion.

\begin{proposition}\label{3.prop}
Let $\iota_k\colon \mathcal{G}^{k}\to \mathrm{Br}(X_k)$,
$f^*_{k+1}$ and $\tilde{f}^*_{k+1}$
be given as above.
Then we have
$\iota_{k+1}\circ f^*_{k+1}=\tilde{f}^*_{k+1}\circ\iota_k$
and that $f^*_{k+1}$ and $\tilde{f}^*_{k+1}$
are order preserving homomorphisms.
\end{proposition}
\begin{proof}
It is not difficult to check,
for every $[g]\in \mathcal{G}^k$,
$$
\left(\iota_{k+1}\circ f^*_{k+1}\right)\left(\left[g\right]\right)
=\left(\tilde{f}^*_{k+1}\circ\iota_k\right)\left(\left[g\right]
\right),
$$
and we have 
$\iota_{k+1}\circ f^*_{k+1}=\tilde{f}^*_{k+1}\circ\iota_k$.

To show that $\tilde{f}^*_{k+1}$ is order preserving,
suppose $[\phi]\in \mathrm{Br}_{\oplus}(X_k)$.
Then there exists a $\psi\in R(X_k)$
such that $\phi\cdot \psi$ is
non-orientation-reversing.
Since $\tilde{f}^*_{k+1}\left(\psi\right)
=\psi\circ f_{k+1}$
is an element of $R(X_{k+1})$ by lemma \ref{3.0}
and $f_{k+1}\colon X_{k+1}\to X_k$
is orientation preserving,
for every orientation preserving parameterized curve
$\gamma \colon \mathbb{R}\to X_{k+1}$,
$f_{k+1}\circ \gamma$ is
an orientation preserving parameterized curve
in $X_k$, and
$$
\left(\left(\phi\circ f_{k+1}\right)\cdot
\left(\psi\circ f_{k+1}\right)\right)
\left(\gamma(t)\right)
=\left(\left(\phi\cdot \psi\right)\circ f_{k+1}
\right)\left(\gamma(t)\right)
=\left(\phi\cdot \psi\right)\circ
\left(f_{k+1}\circ \gamma\right)(t)
$$
does not move in the clockwise direction
as $t\in \mathbb{R}$ increases.
Therefore $\left[\phi\circ f_{k+1}\right]=
\tilde{f}^*_{k+1}\left([\phi]\right)$
is an element of $\mathrm{Br}_{\oplus}(X_{k+1})$,
and $\tilde{f}^*_{k+1}$ is
an order preserving homomorphism.
Since $\iota_k$ is an order preserving isomorphism
by proposition \ref{3.pro},
$f^*_{k+1}={\iota_{k+1}}^{-1}\circ
\tilde{f}^*_{k+1}\circ \iota_k$ 
is also order preserving.
\end{proof}

Then $\{\mathcal{G}^{k},f_{k+1}^*\}$ and
$\{\mathrm{Br}(X_k),\tilde{f}^*_{k+1}\}$ are
directed systems.
Let $\lim\limits_{\longrightarrow} \mathcal{G}^k$
and $\lim\limits_{\longrightarrow} \mathrm{Br}(X_k)$
be  the direct limits of
$\{\mathcal{G}^{k},{f^*_{k+1}}\}$ and
$\{\mathrm{Br}(X_k),\tilde{f}^*_{k+1}\}$,
respectively.

\begin{definition}\label{3.remark}
Recall that $C_+(\mathcal{E}_k,\mathbb{Z})$ is
the subset of $C(\mathcal{E}_k,\mathbb{Z})$
with range in $\mathbb{Z}_+$, and
that $\mathcal{G}^k_+$ is given by
$C_+(\mathcal{E}_k,\mathbb{Z}) / \mathcal{V}_k$.
Since $f^*_{k+1}\colon C(\mathcal{E}_k,\mathbb{Z})\to 
C(\mathcal{E}_{k+1},\mathbb{Z})$ defined by 
$g \mapsto g\circ f_{k+1}$ is an order preserving
homomorphism by proposition \ref{3.prop}, 
$(\lim\limits_{\longrightarrow} \mathcal{G}^k)_+
=\lim\limits_{\longrightarrow} \mathcal{G}_+^k$
is well-defined.
This set, as a positive set, defines the order 
which is the {\it direct limit order}
or the {\it standard order} on
$\lim\limits_{\longrightarrow} \mathcal{G}^k$.
\end{definition}

\subsection*{The standard isomorphism
$\lim\limits_{\longrightarrow} \mathcal{G}^k
\to \mathrm{Br}\left( \overline{X} \right)$.}
Suppose 
$\overline{X}=\lim\limits_{\longleftarrow}X_k$ and
that $\pi_k\colon \overline{X}\to X_k$ is
the projection map to the $k$th coordinate space.
If $\phi$ is an element in $C(X_k, S^1)$,
then $\phi$ induces an element 
$\phi\circ \pi_k \in C(\overline{X},S^1)$,
We will use the isomorphism
$\iota_k\colon \mathcal{G}^k\to \mathrm{Br}(X_k)$
and the natural map
$\mathrm{Br}(X_k)\to \mathrm{Br}(\overline{X})$
defined by $[\phi]\mapsto [\phi\circ \pi_k]$
to make an isomorphism
$\iota\colon
\lim\limits_{\longrightarrow} \mathcal{G}^k
\to \mathrm{Br}\left( \overline{X} \right)$.

Let $1_{X_k}\colon X_k\to S^1$ and
$1_{\overline{X}}\colon \overline{X}\to S^1$
be given by $x_k\mapsto 1\in S^1$ and $x\mapsto 1$
for all $x_k\in X_k$ and $x\in \overline{X}$, respectively.
Suppose that
$\phi$ is an element of $C({X}_k,S^1)$ such that
$\phi$ is homotopic to $1_{X_k}$
by $H\colon X_k\times [0,1]\to S^1$.
Then $\phi\circ \pi_k$ is homotopic to
$1_{\overline{X}}=1_{X_k}\circ\pi_k$
by the map
$\overline{H}\colon \overline{X}\times [0,1]\to S^1$
given by $\overline{H}(x,t)=H(\pi_k(x),t)$.
%
%
%
%
%
%
%
%
%
Thus there is a well-defined map
$$
\pi^*_k\colon 
\mathrm{Br}(X_k)\to \text{Br}(\overline{X})
\text{ given by } [\phi]\mapsto [\phi\circ \pi_k].
$$
Since $(\phi_1\cdot \phi_2)\circ \pi_k
=(\phi_1\circ \pi_k)\cdot (\phi_2\circ \pi_k)$
for all $\phi_1, \phi_2 \in C(X_k,S^1)$,
$\pi^*_k$ is a homomorphism. 
That $f_{k+1}\circ \pi_{k+1}=\pi_k\colon
\overline{X}\to X_k$
implies the following lemma.

\begin{lemma}\label{3.3}
Let $\pi^*_k$ and $\tilde{f}^*_{k+1}$ be defined as above.
Then for all $k$,
$\pi^*_{k+1}\circ \tilde{f}^*_{k+1}=\pi^*_k$.
\end{lemma}

Let $\varphi^*_k\colon \mathrm{Br}(X_k) \to 
\lim\limits_{\longrightarrow} \mathrm{Br}(X_k)$
be the natural map for each $k$. 
If $\varphi^*_k([\phi])=\varphi^*_l([\psi])$
for $[\phi]\in \mathrm{Br}(X_k)$ and
$[\psi]\in \mathrm{Br}(X_l)$,
then there is a positive integer $m\ge k,l$ such that
${\tilde{f}^*_{m+1}\circ\cdots
\circ \tilde{f}^*_{k+1}}([\phi])
={\tilde{f}^*_{m+1}\circ\cdots
\circ \tilde{f}^*_{l+1}}([\psi])$.
Hence 
$$
\pi^*_k([\phi])
=\pi^*_{m+1}\circ \tilde{f}^*_{m+1}\circ\cdots
\circ \tilde{f}^*_{k+1}([\phi])
=\pi^*_{m+1}\circ 
{\tilde{f}^*_{m+1}\circ\cdots
\circ\tilde{f}^*_{l+1}}([\psi])
=\pi^*_l([\psi]),
$$
and there is a well-defined group homomorphism
$$
\pi^*\colon
\lim\limits_{\longrightarrow} \mathrm{Br}(X_k) 
\to \text{Br}(\overline{X})
\text{ given by }
\varphi^*_k([\phi])\mapsto
\pi^*_k([\phi])=\left[\phi\circ \pi_k\right].
$$

\begin{lemma}\label{3.55}
Suppose that $\xi$ is an element of
$C(\overline{X},S^1)$.
Then there exist $\xi^\prime \in C(\overline{X},S^1)$
and $k \ge 0$ such that 
$\xi$  is homotopic to $\xi^\prime$
and $\xi^\prime(x)=\xi^\prime(y)$  if  $x_k=y_k$.
\end{lemma}
\begin{proof}
Define a metric $d$ on $\overline{X}$ by
$
d(x,y)=\sum\limits_{k=0}^{\infty}
\dfrac{1}{2^k}d_k(x_k,y_k)
$
where $x=(x_0,x_1,\dots)$,
$y=(y_0,y_1,\dots)\in \overline{X}$ and
$d_k$ is a  metric on $X_k$ compatible with its topology. 
Since $\overline{X}$ is a compact Hausdorff space,
every element in $C(\overline{X},S^1)$ is
uniformly continuous.
So, for given $\xi$ and $\epsilon > 0$,
there exists a nonnegative integer $k$
such that 
for $x,y\in \overline{X}$,
$x_k=y_k$ implies $d(\xi(x),\xi(y))<\epsilon$.

For $x=(x_0,\dots,x_k,\dots)\in \overline{X}$,
let's denote $x^k=\{ y\in \overline{X}\mid y_k=x_k\}$.
Then for all $a,b\in x^k$,
$d(\xi(a),\xi(b))<\epsilon$
and we can choose a point $\tilde{x}\in S^1$ 
such that $\tilde{x}$ is the center of 
the smallest interval
containing $\xi(x^k)$ in $S^1$.
Define $\xi^\prime\colon \overline{X}\to S^1$ by
$\xi^\prime|_{x^k}=\tilde{x}$.
Then it is clear that
$\xi^\prime\in C(\overline{X},S^1)$ and 
$\xi^\prime(x)=\xi^\prime(y)$ if $x_k=y_k$.
Since $d(\xi(x),\xi^\prime(x))<\epsilon$
for all $x\in \overline{X}$,
$\xi$ is homotopic to $\xi^\prime$.
\end{proof}

\begin{proposition}\label{3.4}
Let $\pi^*$ be defined as above.
Then $\pi^*$ is a group isomorphism.
\end{proposition}
\begin{proof}
To show that $\pi^*$ is surjective,
suppose $\xi\in C(\overline{X}, S^1)$ 
and that $\xi^\prime$ and $k$ are given in
Lemma \ref{3.55}.
Define $\phi_k\colon X_k\to S^1$ by
$x_k\mapsto \xi^\prime(x)$
for $x=(x_0,\dots,x_k,\dots)\in \overline{X}$.
Then $\phi_k$ is well-defined,
and it is trivial that 
$\phi_k\circ \pi_k=\xi^\prime$.
Therefore $\xi\in C(\overline{X},S^1)$
is homotopic to $\phi_k\circ \pi_k$,
and $\pi^*\colon 
\lim\limits_{\longrightarrow} \mathrm{Br}(X_k) 
\to \text{Br}(\overline{X})$
is surjective.

Suppose  $\xi_1, \xi_2 \in C(\overline{X},S^1)$
and that $\xi_1$ is homotopic to $\xi_2$.
Then by the surjectivity of $\pi^*$,
there exist nonnegative integers $k \le l$
and $\phi\in C(X_k,S^1)$, $\psi\in C(X_l,S^1)$
such that $\xi_1$ is homotopic to $\phi\circ \pi_k$
and $\xi_2$ is homotopic to $\psi\circ \pi_l$.
Since $\phi\circ \pi_k
=\phi\circ f_{k+1}\circ\cdots\circ f_l\circ\pi_l$,
we have
$$
\varphi_l^*\left([\psi]\right)=
\varphi_l^*\left(
[\phi\circ f_{k+1}\circ\cdots\circ f_l]
\right)
=\varphi_l^*\circ \tilde{f}_l^*\circ\cdots
\circ \tilde{f}^*_{k+1}
\left( [\phi] \right)
=\varphi_k^*\left( [\phi] \right).
$$
Hence $\pi^*$ is injective.
\end{proof}

Therefore the isomorphisms
$\iota_k\colon \mathcal{G}^K\to \mathrm{Br}(X_k)$
and $\pi^*\colon 
\lim\limits_{\longrightarrow} \mathrm{Br}(X_k) 
\to \text{Br}(\overline{X})$
induce an isomorphism
$\iota\colon
\lim\limits_{\longrightarrow} \mathcal{G}^k
\to \text{Br}(\overline{X})$.

\subsection*{Order isomorphism}
Assume now that the presentation $\{X_k,f_k \}$ satisfies 
the following Simplicity Condition:
\begin{itemize}
\item[ ]
for each $k\ge 1$ there exists $\kappa(k)\ge k$ such that
for every $l\ge \kappa(k)$ and $e\in \mathcal{E}_l$
$f_{k+1}\circ\cdots\circ f_l(e)=X_k$
where $\mathcal{E}_l$ is the edge set of $X_l$.
\end{itemize}
Then the winding order on
$\mathrm{Br}(X_k)$ and $\mathrm{Br}(\overline{X})$
is an order.

\begin{theorem}\label{3.6}
Suppose that the presentation $\{X_k,f_k \}$ satisfies 
the above Simplicity Condition.
Then $\iota\colon
\left(
\lim\limits_{\longrightarrow} \mathcal{G}^k,
\lim\limits_{\longrightarrow} \mathcal{G}^k_+\right)
\to 
\left( \text{\rm Br}(\overline{X}),
\text{\rm  Br}_{\oplus}(\overline{X})\right)$
is an isomorphism of ordered groups.
\end{theorem}
\begin{proof}
(\emph{Trivial case}).
Suppose that 
all but finitely many $X_k$ has a unique edge,
i.e., 
$X_k$ is homeomorphic to 
a circle $S^1$ with a unique vertex
by the standing assumption \ref{2sa},
and that the connection map $f_k\colon X_k\to X_{k-1}$
is the identity map if $X_k=X_{k-1}=S^1$.
Then it is obvious that
$$
\left(
\lim\limits_{\longrightarrow} \mathcal{G}^k,
\lim\limits_{\longrightarrow} \mathcal{G}^k_+\right)
\cong
\left( \text{\rm Br}(\overline{X}),
\text{\rm  Br}_{\oplus}(\overline{X})\right)
=\left( \mathbb{Z},\mathbb{Z}_+\right),
$$
and $\iota$ is an isomorphism.

\noindent
(\emph{Nontrivial case}).
We have that $\iota$ is a group isomorphism, and clearly
$\iota(\lim\limits_{\longrightarrow} \mathcal{G}^k_+)
\subseteq \text{Br}_{\oplus}(\overline{X})$.
It remains to show that $\iota$ maps
$\lim\limits_{\longrightarrow} \mathcal{G}^k_+$ onto
$\text{Br}_{\oplus}(\overline{X})$. 
So we assume that $[\phi]$ is an element of
$\text{Br}_{\oplus}(\overline{X})$.
Then there is an $[h]$ in $\mathcal{G}^k$
for some $k\ge 0$ such that
$[\phi]=[\phi_h\circ\pi_k]$,
and we need to show $[h] \in \mathcal{G}^k_+$.

That $[\phi]$ is an element of
$\text{Br}_{\oplus}(\overline{X})$ 
implies that
there is a map $\gamma\in R(\overline{X})$ such that
$(\phi_h\circ\pi_k)\cdot\gamma$ is
non-orientation-reversing. 
Since $\gamma$ is an element of $R(\overline{X})$,
there is a continuous map $g\colon \overline{X}\to \mathbb{R}$
such that $\gamma(x)=\exp(2\pi i g(x))$.
For $y=(y_0,\dots,y_k,\dots)\in \overline{X}$,
if $y_k=e(t)$ for $e\in \mathcal{E}_k$ and $t\in[0,1]$,
then $\phi_h\circ\pi_k\cdot \gamma$ is defined by
$y \mapsto \exp(2\pi i (t h(e)+g(y)) )$,

Suppose that $(\phi_h\circ\pi_k)\cdot\gamma$
is a constant map to $S^1$.
Then we have 
$[(\phi_h\circ\pi_k)\cdot\gamma]=
[\phi_h\circ\pi_k]\cdot[\gamma]=[\phi_h\circ\pi_k]=[1]$
in $\text{Br}(\overline{X})$ 
as $\gamma$ is homotopic to the identity element in
$\text{Br}(\overline{X})$.
Hence  the equivalence class of $h$ is
the identity element in
$\lim\limits_{\longrightarrow} \mathcal{G}^k$, for
$\iota \colon \lim\limits_{\longrightarrow} \mathcal{G}^k
\to \text{Br}(\overline{X})$ is an isomorphism.

Next suppose that $(\phi_h\circ\pi_k)\cdot\gamma$
is not constant on $S^1$.
Then there are nonnegative integer $m$,
a small 
interval $I$
contained in some edge $e^\prime$ of $X_{k+m}$,
and $\epsilon >0$ such that
if  $\Gamma$ is any orientation preserving curve
in $\overline{X}$  and $\pi_{k+m}(\Gamma|_{[a,b]}) = I$,
then
$length\, 
\{((\phi_h\circ\pi_k)\cdot\gamma)\circ \Gamma|_{[a,b]}\} 
>\epsilon$.

Given an arbitrary constant $L$, 
by the simplicity condition
we can choose a sufficiently large integer $M$
such that $e^\prime$ is covered under
$f^{k+m+1}\circ\cdots\circ f^{k+m+M}$
at least $L$ times by every edge in $\mathcal{E}_{k+m+M}$.

Define 
$H=f^*_{k+m+M}\circ\cdots \circ f^*_{k+1}(h)=
h\circ f_{k+1}\circ\cdots\circ f_{k+m+M}
\in C(\mathcal{E}_{k+m+M},\mathbb{Z})$.
Then by Lemma~\ref{3.3},
$\phi_H\circ\pi_{k+m+M}\in C(\overline{X},S^1)$ 
is homotopic to $\phi_h\circ\pi_k$.
For $x=(x_0,\dots,x_{k+m+M},\dots)\in \overline{X}$,
as $x_{k+m+M}$ moves forward through
a directed edge $e$ of $\mathcal{E}_{k+m+M}$,
its image under $\phi_H\circ\pi_{k+m+M}$ moves 
$\sum h(\hat e)\cdot n_e(\hat e)$ times
around $S^1$ where $n_e(\hat e)$ is the number of times  
$e$ covers $\hat e \in \mathcal{E}_k$ under the map
$f^{k+1}\circ\cdots\circ f^{k+m+M}$.

\begin{lemma}\label{3.7} 
For every edge $e\in  \mathcal{E}_{k+m+M}$,
$H(e)\ge 2\pi L\epsilon - 2\max |g|$.
\end{lemma}
\begin{proof}[Proof of Lemma]
Regard $e$ as a curve $e(t)$, $0\leq t\leq 1$,
and pick a curve
$\Gamma\colon [0,1]\to \overline{X}$ 
such that $\pi_{k+m+M}\circ \Gamma (t) = e(t)$. 
As $t$ increases from $0$ to $1$, 
the point 
$$
\big( (\phi_h\circ \pi_k) \cdot \gamma \big) \circ \Gamma (t) 
=
\big( \phi_h \circ \pi_k \circ \Gamma(t) \big) \cdot 
\big( \gamma \circ \Gamma (t) \big) 
$$
moves counterclockwise on $S^1$ from 
$e^{2\pi ig(\Gamma (0))}$ to 
$e^{2\pi i\left(g(\Gamma (1))+ H(e)\right)}$, 
covering an arclength $A$ in the plane 
such that 
$$
A \le \frac{1}{2\pi}\left( H(e) + 2\text{max} |g|\right).
$$
Because 
$\phi_h \circ \pi_k \circ \Gamma =
\phi_h \circ f_{k+1} \circ \cdots \circ f_{k+m+M} 
\circ \pi_{k+m+M} \circ \Gamma$, 
as $t$ runs from $0$ to $1$ the curve 
$f_{k+m+1}\circ \cdots \circ f_{k+m+M}\circ
\pi_{k+m+M}\circ \Gamma (t)$ 
wraps around $e^\prime$ at least $L$ times,
and therefore $A\geq L\epsilon$. Consequently 
$2\pi L\epsilon - 2{\max}|g| \leq H(e)$ as required. 
\end{proof}

Since we can choose $M$
to make $L$ as large as we wish,
we can make the choice to guarantee
$H(e)>0$ for every edge.
Therefore
$[H]=[h]$ is an element of $\mathcal{G}^k_{+}$.
\end{proof}

\subsection*{Dimension group}
Let $M$ be an $r\times s$ nonnegative integer matrix.
Then the matrix $M$ determines a homomorphism
$\mathbb{Z}^s\to \mathbb{Z}^r$
by the ordinary matrix multiplication.
The {\it simplicial order} on $\mathbb{Z}^r$
is the usual ordering
$\mathbb{Z}^r_+=\{(n_1,\dots,n_r)\mid n_i\ge 0\}$.
Then the corresponding homomorphism
$M\colon \mathbb{Z}^s\to \mathbb{Z}^r$
is {\it positive} with respect to the simplicial order,
that is, $a\ge 0$ implies $M(a)\ge 0$.

\begin{definition}[{\cite[\S2]{ef}}]\label{3.ef}
Let $M_i$ be an $r(i)\times r(i-1)$
nonnegative integer matrix.
For a system of ordered groups and positive maps
$$
\mathbb{Z}^{r(0)}\overset{M_1}{\longrightarrow}
\mathbb{Z}^{r(1)}\overset{M_2}{\longrightarrow}\cdots ,
$$
the set theoretic direct limit
$\lim\limits_{\longrightarrow}
\left(\mathbb{Z}^{r(i)},M_i\right)$
is an ordered group under the usual limit addition operation
with the positive cone
$\lim\limits_{\longrightarrow}
\left(\mathbb{Z}^{r(i)}_+,M_i\right)
=\bigcup\limits_{i=1}^{\infty}M_{i\infty}
\left( \mathbb{Z}^{r(i-1)}_+ \right)$
where $M_{i\infty}$ is the induced map
from $\mathbb{Z}^{r(i-1)}$ to
the direct limit
$\lim\limits_{\longrightarrow}
\left(\mathbb{Z}^{r(i)},M_i\right)$.

An ordered group $(G,G_+)$ is called
a {\bf dimension group}
if it is order isomorphic to
the limit of a system of simplicially ordered groups
with positive maps.

Let $(G,G_+)$ be a dimension group.
A subgroup $H$ of $G$ is
called an {\it order ideal}
if $H$ is an ordered group with the positive cone
$H_+=H\cap G_+$
and $0\le a \le b \in H$ implies $a\in H$.
The dimension group $(G,G_+)$ is called {\it simple}
if it has no proper order ideal.

In a simple dimension group $(G,G_+)$ with 
an element $g\in G$,
if neither $g$ nor $-g$ lies in $G_+$,
then $g$ is called an {\it infinitesimal} element.
If $u$ is an order unit  and $g$ is
an infinitesimal element of $G$,
then $g+u$ is also an order unit. 
\end{definition}

It is well known that
a dimension group defined as above
by matrices $M_i$ is simple
if for every $i$ there exists $j$ such that
all entries of the matrix
$M_jM_{j-1}\cdots M_{i+1}M_i$ are strictly positive.

Suppose that $\{X_k,f_k\}$ is a presentation of
an (orientable) branched matchbox manifold
with the edge set $\mathcal{E}_k$ of $X_k$.
Then for each edge $e_i\in \mathcal{E}_k$,
$f_k(e_i)$ is a path 
$e_{i,1}^{s(1)}\cdots e_{i,j(i)}^{s(j(i))}$ in $X_{k-1}$
such that
$s(j)=\pm 1$ denotes the direction and
the terminal point of $e_{i,j}^{s(j)}$ is 
the initial point of $e_{i,j+1}^{s(j+1)}$
for $1\le j<j(i)$.
Therefore we can define an induced map
$\check{f}_k\colon \mathcal{E}_k\to \mathcal{E}_{k-1}^*$
by
$$
\check{f}_k\colon e_i\mapsto
e_{i,1}^{s(1)}\cdots e_{i,j(i)}^{s(j(i))}.
$$

\begin{definition}
Suppose that $X_k$ has $n_k$ edges for all $k\ge 0$.
Then the {\bf adjacency matrix} $M_k$ of
$\left(
\check{f}_k,\mathcal{E}_k,\mathcal{E}_{k-1}
\right)$
is an $n_k\times n_{k-1}$ matrix such that
for any edges $e_i\in \mathcal{E}_k$ and 
$e_j\in \mathcal{E}_{k-1}$,
$M_k(i,j)$ is the number of times $\check{f}_k(e_i)$
covers $e_j$ ignoring the direction of the covering.
\end{definition}

\begin{lemma}[{\cite[\S3]{ef}}]\label{3.h}
A countable  ordered group is a dimension group
if and only if
it is unperforated and
has the Riesz Interpolation Property.
\end{lemma}
 

\begin{proposition}\label{3.88}
Suppose that $\{X_k,f_k\}$ is a presentation
of a compact connected orientable branched matchbox manifold
with the adjacency matrices $M_k$.
Then 
\begin{itemize}
\item[(1)]
$\left(
\lim\limits_{\longrightarrow}
C(\mathcal{E}_k,\mathbb{Z}),
\lim\limits_{\longrightarrow}
C_+(\mathcal{E}_k,\mathbb{Z})
\right)
\cong
\left(
\lim\limits_{\longrightarrow}
\left(\mathbb{Z}^{n_k},M_k\right),
\lim\limits_{\longrightarrow}
\left(\mathbb{Z}^{n_k}_+,M_k\right)
\right)$.
\end{itemize}
If the presentation 
satisfies the Simplicity Condition, then 
\begin{itemize}
\item[(2)]
$\left(\lim\limits_{\longrightarrow}
C(\mathcal{E}_k,\mathbb{Z}),
\lim\limits_{\longrightarrow}
C_+(\mathcal{E}_k,\mathbb{Z})
\right)$ and
$\left(\lim\limits_{\longrightarrow} \mathcal{G}^k,
\lim\limits_{\longrightarrow} \mathcal{G}^k_+ \right)$
are simple dimension groups.
\end{itemize}
\end{proposition}
\begin{proof}
(1)
For each $g\in C(\mathcal{E}_{k-1},\mathbb{Z})$
and $f^*_k\colon C(\mathcal{E}_{k-1},\mathbb{Z})
\to C(\mathcal{E}_k,\mathbb{Z})$ given by 
$g\mapsto g\circ f_k$,
if we represent $g$ as
$\left(g(e_1),\dots , g(e_{n_{k-1}})\right)
\in \mathbb{Z}^{n_{k-1}}$,
then $C(\mathcal{E}_{k-1},\mathbb{Z})$
is isomorphic to $\mathbb{Z}^{n_{k-1}}$
and 
$f^*_k(g)=g\circ f_k$ is given by
$M_k\cdot \left(
g(e_1),\dots , g(e_{n_{k-1}})
\right)^t$.
Hence we have
$\lim\limits_{\longrightarrow}
C(\mathcal{E}_{k-1},\mathbb{Z})
\cong
\lim\limits_{\longrightarrow}
\left(\mathbb{Z}^{n_k}, M_k\right)$.
Since $C_+(\mathcal{E}_{k-1},\mathbb{Z})$
is the set of elements in
$C(\mathcal{E}_{k-1},\mathbb{Z})$
with range in $\mathbb{Z}_+$,
$C(\mathcal{E}_{k-1},\mathbb{Z})$
is simplicially ordered,
and so is
$\lim\limits_{\longrightarrow}
C(\mathcal{E}_{k},\mathbb{Z})$.
Therefore
$\left(\lim\limits_{\longrightarrow}
C(\mathcal{E}_{k},\mathbb{Z}),
\lim\limits_{\longrightarrow}
C_+(\mathcal{E}_{k},\mathbb{Z})
\right)$
is order isomorphic to
$\left(
\lim\limits_{\longrightarrow}
\left(\mathbb{Z}^{n_k},M_k\right),
\lim\limits_{\longrightarrow}
\left(\mathbb{Z}^{n_k}_+,M_k\right)
\right)$.

%
%

\noindent
(2) 
Suppose that $H$ is a proper order ideal
of $\left(\lim\limits_{\longrightarrow}
C(\mathcal{E}_{k},\mathbb{Z}),
\lim\limits_{\longrightarrow}
C_+(\mathcal{E}_{k},\mathbb{Z}) \right)$
and that $b\in H_+$.
Then there exist a nonnegative integer $k$ and
$h\in C_+(\mathcal{E}_k,\mathbb{Z})$
such that $b=\left[ h \right] \in
\lim\limits_{\longrightarrow}
C(\mathcal{E}_{k},\mathbb{Z})$.
By the Simplicity Condition,
there is a nonnegative integer $\kappa(k)\ge k$
such that, for every $l\ge \kappa(k)$ and
every edge $e\in \mathcal{E}_l$,
$f_{k+1}\circ \cdots \circ f_l(e)=X_k$.
If $a\in\lim\limits_{\longrightarrow}
C_+(\mathcal{E}_{k},\mathbb{Z}) $,
then we can choose a positive integer $l\ge \kappa(k)$
and $g\in C_+(\mathcal{E}_l,\mathbb{Z})$
such that $a=\left[ g \right]$.
Let $n=\max\limits_{e\in \mathcal{E}_l} g(e)$.
Then $n\cdot b 
=\left[n\cdot f^*_{l}\circ\cdots\circ f^*_{k+1}\circ h\right]
\in H_+$ and
$n\cdot f^*_{l}\circ\cdots\circ f^*_{k+1}\circ h
-g\in C_+(\mathcal{E}_l,\mathbb{Z})$.
So we have $0\le a\le n\cdot b$ and $a\in H_+$. 
Therefore $H_+
=\lim\limits_{\longrightarrow}
C_+(\mathcal{E}_{k},\mathbb{Z})$,
and $\left(
\lim\limits_{\longrightarrow}C(\mathcal{E}_{k},\mathbb{Z}),
\lim\limits_{\longrightarrow}C_+(\mathcal{E}_{k},\mathbb{Z}) 
\right)$
is a simple dimension group.

The group $\left(\lim\limits_{\longrightarrow} \mathcal{G}^k,
\lim\limits_{\longrightarrow} \mathcal{G}^k_+ \right)
\cong
\left( \mathrm{Br}(\overline{X}),
\mathrm{Br}_{\oplus}(\overline{X})\right)$
is an unperforated ordered group by proposition \ref{2.r},
and its positive set is the image of the positive set of 
$\lim\limits_{\longrightarrow} C(\mathcal{E}_{k},\mathbb{Z})$
under the quotient map
$\chi\colon
\lim\limits_{\longrightarrow}
C(\mathcal{E}_{k},\mathbb{Z}) \to
\lim\limits_{\longrightarrow} \mathcal{G}^k$.
We claim that with this quotient order,
$\left(\lim\limits_{\longrightarrow} \mathcal{G}^k,
\lim\limits_{\longrightarrow} \mathcal{G}^k_+ \right)$
satisfies the Riesz Interpolation Property
(and therefore by Lemma \ref{3.h} is a dimension group.)
(We learned this argument from unpublished remarks
of David Handelman.
The general line of argument is also
implicit in remarks on pp.\:58 and 66 of \cite{gps}.)



Let $V=\text{ker}\,\chi$.
Note that
if $V$ contains a nonzero positive element $u$,
then for every $g\in \lim\limits_{\longrightarrow}
C_+(\mathcal{E}_{k},\mathbb{Z})$
we have for some integer $n$ that
$0\le g \le nu$,
and therefore $0\le \chi(g) \le 0$,
which contradicts the image of $\chi$ being
a nontrivial ordered group.
Therefore all elements of $V$ are infinitesimals.

To show the Riesz Interpolation Property,
suppose that $[a_1],[a_2],[b_1],[b_2]\in 
\lim\limits_{\longrightarrow} \mathcal{G}^k$
satisfy $[a_i] < [b_j]$ ($i,j=1,2$).
Let $a_i$ and $b_j\in \lim\limits_{\longrightarrow}
C(\mathcal{E}_{k},\mathbb{Z})$
be preimages of $[a_i]$ and $[b_j]$, respectively.
Since $-[a_i]+[b_j]$ is a nonzero positive element of
$\lim\limits_{\longrightarrow} \mathcal{G}^k$,
there exists a $v_{i,j}\in V$ such that
$-a_i+v_{i,j}+ b_j$ is a nonzero positive element of
$\lim\limits_{\longrightarrow}
C(\mathcal{E}_{k},\mathbb{Z})$.
Because $v_{ij}$ is an infinitesimal element,
it follows that
$-a_i+ b_j$ is a nonzero positive element of
$\lim\limits_{\longrightarrow}C(\mathcal{E}_{k},\mathbb{Z})$,
and $a_i < b_j$ for $i,j=1,2$.
Hence by the Riesz Interpolation Property
for $\lim\limits_{\longrightarrow}
C(\mathcal{E}_{k},\mathbb{Z})$
there exists an element $c\in \lim\limits_{\longrightarrow}
C(\mathcal{E}_{k},\mathbb{Z})$
such that
$a_i\leq c \leq b_j$.
Then by definition of the quotient order
we have $[a_i] \leq [c] \leq [b_j]$
for all $i,j$ as required. 
Therefore
$\left(\lim\limits_{\longrightarrow} \mathcal{G}^k,
\lim\limits_{\longrightarrow} \mathcal{G}^k_+ \right)$
is a dimension group by lemma \ref{3.h}.

Suppose that $\left(G,G_+\right)$ is a proper order ideal
of $\left(\lim\limits_{\longrightarrow} \mathcal{G}^k,
\lim\limits_{\longrightarrow} \mathcal{G}^k_+ \right)$.
Then it is not difficult to see that
$\left(H,H_+\right)=\left(\chi^{-1}(G),\chi^{-1}(G_+)\right)$
is a proper order ideal of 
$\left(
\lim\limits_{\longrightarrow}C(\mathcal{E}_{k},\mathbb{Z}),
\lim\limits_{\longrightarrow}C_+(\mathcal{E}_{k},\mathbb{Z}) 
\right)$ which is a simple dimension group.
Therefore  $\left(\lim\limits_{\longrightarrow} \mathcal{G}^k,
\lim\limits_{\longrightarrow} \mathcal{G}^k_+ \right)$
is a simple dimension group.
\end{proof}

If each graph $X_k$ is a wedge of circle,
then $V_k=\{ 0 \}$ as each edge in $X_k$
is a cycle.
So we have the following corollary

\begin{corollary}\label{3.8888}
Suppose that the presentation
$\{X_k,f_k\}$ satisfies the Simplicity Condition
and that each graph $X_k$ is a wedge of circles.
Then $\left(\lim\limits_{\longrightarrow} \mathcal{G}^k,
\lim\limits_{\longrightarrow} \mathcal{G}^k_+ \right)$
is order isomorphic to
$\left(\lim\limits_{\longrightarrow}
\left(\mathbb{Z}^{n_k},M_k\right),
\lim\limits_{\longrightarrow}
\left(\mathbb{Z}^{n_k}_+,M_k\right)\right)$.
\end{corollary}

The following corollary follows from
Observation \ref{1.1} and Theorem \ref{3.6}.

\begin{corollary}\label{3.888}
Suppose that $(\overline{X_i},\overline{f_i})$
is a compact connected orientable branched matchbox manifold 
with the Simplicity Condition
for $i=1,2$. 
If $\overline{X_1}$ is homeomorphic to $\overline{X_2}$,
then $\lim\limits_{\longrightarrow} \mathcal{G}^k_1$
is order isomorphic to
$\lim\limits_{\longrightarrow} \mathcal{G}^k_2$.
\end{corollary}

\begin{remark}
\begin{itemize}
\item[(1)]
The dimension group of adjacency matrices
is not a homeomorphism invariant.
See example \ref{4.ex}.
\item[(2)]
The isomorphism in corollary \ref{3.888}
need not respect distinguished order unit
(\cite[\S1]{bh}).
\end{itemize}
\end{remark}

\section{One-dimensional generalized solenoid}
One interesting class of branched matchbox manifolds is
one-dimensional branched solenoids,
including one-dimensional
generalized solenoids of Williams
(\cite{w1, w2, yi}).
Let $X$ be a directed graph with vertex set
$\mathcal{V}$ and edge set $\mathcal{E}$,
and $f\colon X\to X$ a continuous map.
We define some axioms which
might be satisfied by $(X,f)$ (\cite{yi}).
\begin{enumerate}
\item[Axiom 0.]({\it Indecomposability})
$(X,f)$ is indecomposable.
\item[Axiom 1.]({\it Nonwandering})
All points of $X$ are nonwandering under $f$.
\item[Axiom 2.] ({\it Flattening}) 
There is $k\ge 1$ such that  for all $x\in X$ 
there is an open neighborhood $U$ of $x$
such that $f^k(U)$ is homeomorphic to
$(-\epsilon, \epsilon )$.
\item[Axiom 3.]({\it Expansion}) 
There are a metric $d$ compatible with the topology
and positive constants $C$ and $\lambda$ with 
$\lambda >1$ such that for all $n>0$ and 
all points $x,y$ on a common edge of $X$, 
if $f^n$ maps the interval $[x,y]$ into an edge,
then $d(f^nx,f^ny)\geq C\lambda^n d(x,y)$.  
\item[Axiom 4.] ({\it Nonfolding})
$f^n|_{X-\mathcal{V}}$ is locally one-to-one
for every positive integer $n$.
\item[Axiom 5.]({\it Markov})
$f(\mathcal{V})\subseteq \mathcal{V}$.
\end{enumerate}

Let $\overline{X}$ be the inverse limit space
$$
\overline{X}
=X\overset{f}{\longleftarrow}
X\overset{f}{\longleftarrow}\cdots
=\bigl\{(x_0,x_1,x_2,\dots )\in
       \prod_0^\infty X \, |\, f(x_{n+1})=x_n \bigr\},
$$
and $\overline{f}\colon \overline{X}\to \overline{X}$
the induced homeomorphism defined by
$$
(x_0,x_1,x_2,\dots )\mapsto 
(f(x_0),f(x_1),f(x_2),\dots )=(f(x_0),x_0,x_1,\dots).
$$

Let $Y$ be a topological space and 
$g\colon Y\to Y$ a homeomorphism.
We call $Y$ a {\bf 1-dimensional generalized solenoid}
or {\bf $1$-solenoid} and $g$ a {\bf solenoid map}
if there exist a directed graph $X$
and a continuous map $f\colon X\to X$
such that $(X,f)$ satisfies all six Axioms
and $(\overline{X},\overline{f})$
is topologically conjugate to $(Y,g)$.
If $(X,f)$ satisfies all Axioms
except possibly the Flattening Axiom,
then we call $Y$ a {\bf branched solenoid}. 
If we can choose
the direction of each edge in $X$ so that
the connection map $f\colon X\to X$ is
orientation preserving, then we call
$({X},{f})$ an {\it orientable presentation},
and $Y$ an {\it orientable} (branched) solenoid.
If $(Y,g)$ is a branched solenoid
with a presentation $({X},{f})$,
then there exist an $n\times n$
adjacency matrix $M_{X,f}$
where $n$ is the cardinal number of
the set of edges in $X$.
If $X$ is a wedge of circles and 
$f$ leaves the unique branch point of $X$ fixed,
then we say $({X},{f})$ an
{\it elementary presentation}.





We have the following proposition
from theorem \ref{3.6} and
corollary \ref{3.8888}.

\begin{proposition}
Suppose that $\left(\overline{X},\overline{f}\right)$ is
an orientable branched solenoid
with an adjacency matrix $M$.
Then $\iota\colon
\left(\lim\limits_{\longrightarrow}
\left(\mathbb{Z}^n, M\right),
\lim\limits_{\longrightarrow}
\left(\mathbb{Z}^n_+, M\right)\right)
\to 
\left( \text{\rm Br}(\overline{X}),
\text{\rm  Br}_{\oplus}(\overline{X})\right)$
is an epimorphism of ordered groups.
If $(X,f)$ is an elementary presentation,
then $\iota$ is an isomorphism.
\end{proposition}

\begin{remark}
We need the elementary presentation condition for
the injectivity of $\iota$.
See example \ref{4.ex}.
\end{remark}

\begin{example}[{\cite[\S2]{yi}} and {\cite[\S7.5]{lm}}]
Let $X$ be the unit circle on the complex plane.
Suppose that $1$ and $-1$ are the vertices of $X$,
and that the upper half circle $e_1$
and the lower half circle $e_2$
with counterclockwise direction are the edges of $X$.
Define $f\colon X\to X$ by $f\colon z\mapsto z^2$.
The $\check{f}\colon \mathcal{E}_X\to \mathcal{E}_X^*$
is given by
$\check{f}\colon
e_1\mapsto e_1e_2,\,\, e_2\mapsto e_1e_2$,
and the adjacency matrix is
$$
M_{X,f}=
\begin{pmatrix}
1&1\\1&1
\end{pmatrix}.
$$
Therefore we have
$$
\left( \text{\rm Br}(\overline{X}),
\text{\rm Br}_{\oplus}(\overline{X})\right)
=\left(\mathbb{Z}\left[1/2\right],
\mathbb{Z}\left[1/2\right]\cap \mathbb{R}_+\right).
$$

We give Figure \ref{xcv1}
to represent
the presentation $(X,f)$ with the wrapping rule $\tilde{f}$.
\begin{figure}[ht]
\centerline{\includegraphics[height=2.9cm]{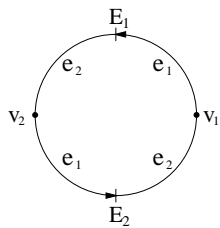}}
\caption{$(X,f)$ with the wrapping rule $\tilde{f}$.}\label{xcv1}
\end{figure}

Similarly,
if $(Y,g)$ is given by Figure \ref{or6},
\begin{figure}[h]
\centerline{\scalebox{.28}{\includegraphics{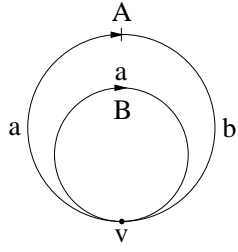}}}
\caption{$(Y,g)$ with wrapping rule $\check{g}$.}\label{or6}
\end{figure}
then $(Y,g)$ does not satisfy the flattening axiom
and $(\overline{Y},\overline{g})$ is a branched solenoid.
The wrapping rule 
$\check{g}\colon \mathcal{E}_Y\to \mathcal{E}_Y^*$
is given by $a\mapsto ab$, $b\mapsto a$
and the adjacency matrix is
$$
M_{Y,g}=
\begin{pmatrix}
1&1\\1&0
\end{pmatrix}.
$$
Thus $\text{\rm Br}(\overline{Y})=\mathbb{Z}\oplus \mathbb{Z}$
and
$\text{\rm Br}_{\oplus}(\overline{Y})=
\left\{
\textbf{v} \in \mathbb{Z}\oplus \mathbb{Z} \mid 
\textbf{v} \cdot \left(\frac{1+\sqrt{5}}{2},1\right) > 0
\right\}
\cup \{ \textbf{0} \}$.
\end{example}

The following example shows that
the dimension group of adjacency matrices
induced by a presentation is not
a homeomorphism invariant.

\begin{example}[{\cite[4.8 and 5.1]{yi}}]\label{4.ex}
Let $X$ be a wedge of two circles $a,b$
with a unique vertex $p$, and
$f\colon X\to X$ defined by
$a\mapsto aab$ and $b\mapsto ab$.
So $(X,f)$ is given by 
Figure~\ref{xcv3}.
\begin{figure}[h]
\centerline{\scalebox{.195}{\includegraphics{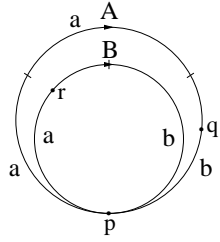}}}
\caption{$(X,f)$ with a unique vertex $\{p\}$}\label{xcv3}
\end{figure}
Suppose that $Y$ is given by 
Figure \ref{xcv5}
\begin{figure}[h]
\centerline{\scalebox{.35}{\includegraphics{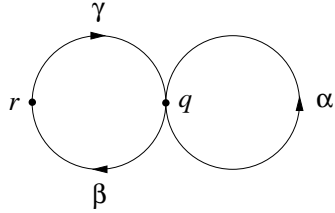}}}
\caption{The graph $Y$ with two vertices $\{q,r\}$}
\label{xcv5} 
\end{figure}
and that the wrapping rule
$\check{g}\colon \mathcal{E}_Y\to \mathcal{E}_Y^*$
is given by
$$
\alpha \mapsto \gamma\alpha\beta,\quad
\beta \mapsto \gamma, \quad
\gamma \mapsto \beta\gamma\alpha\beta.
$$
Then it is shown in \cite[4.8]{yi}
that $\left(\overline{X},\overline{f}\right)$
is topologically conjugate to
$\left(\overline{Y},\overline{g}\right)$.
And their adjacency matrices are
given by the following  matrices
\[
M_{(X,f)}=
\begin{pmatrix}
2&1\\ 
1&1
\end{pmatrix}
\text{ and }
M_{(Y,g)}=
\begin{pmatrix}
1&1&1\\ 
0&0&1\\
1&2&1
\end{pmatrix}.
\]
Since the determinants of $M_{(X,f)}$
and $M_{(Y,g)}$ are $1$ and $-1$, respectively,
$M_{(X,f)}$ and $M_{(Y,g)}$ are
invertible over $\mathbb{Z}$.
Hence the dimension group of $M_{(X,f)}$ is $\mathbb{Z}^2$
and that of $M_{(Y,g)}$ is $\mathbb{Z}^3$.
Therefore the dimension group of $M_{(X,f)}$
is not isomorphic to the dimension group of
$M_{(Y,g)}$

Since $(X,f)$ is elementary presented,
the dimension group of $M_{(X,f)}$ is
order isomorphic to the Bruschlinsky group of
$(\overline{X},\overline{f})$.
And 
the Bruschlinsky group of
$(\overline{Y},\overline{g})$ is
given by the dimension group of
$\begin{pmatrix} 1&1\\1&2 \end{pmatrix}$.
Hence we have
$\mathrm{Br}(\overline{X}) \cong
\mathrm{Br}(\overline{Y}) \cong
\mathbb{Z}\oplus \mathbb{Z}$
with
$$
\mathrm{Br}_{\oplus}(\overline{X})\cong
\mathrm{Br}_{\oplus}(\overline{Y})\cong
\left\{
\textbf{v} \in \mathbb{Z}\oplus \mathbb{Z} \mid 
\textbf{v} \cdot \left(\frac{1+\sqrt{5}}{2},1\right) > 0
\right\}
\cup \{ \textbf{0} \}.
$$
\end{example}

\begin{ack}
This paper is a part of my Ph.D. dissertation
at UMCP.
I sincerely express my gratitude 
to my advisor Dr.{\;}M.{\;}Boyle
for his encouragement and advice.
The proof of Lemma~\ref{3.7} was suggested
by Dr.{\;}Boyle.
By kind permission, I presented his proof here.
\end{ack}

\bibliographystyle{amsplain}

\end{document}